\documentclass{article}

\usepackage{amsfonts}
\usepackage{amssymb}
\usepackage{amsmath}
\usepackage[mathcal]{euscript}
\usepackage[all]{xypic}

\newcommand{\CC}{\mathbb{C}}
\newcommand{\CP}{\mathbb{CP}}
\newcommand{\QQ}{\mathbb{Q}}

\newcommand{\NN}{\mathbb{N}}
\newcommand{\ZZ}{\mathbb{Z}}
\newcommand{\sph}{\mathbb{S}}

\newcommand{\bq}{/\!/}
\newcommand{\smsh}{\wedge}
\newcommand{\wsm}{\vee}

\title{Complex Cobordism vs Representing Formal Group Laws}
\author{Jesse C. McKeown\footnote{Interrupted in his doctoral studies at University of Toronto.}}
\date{12 May, 2016}

\begin{document}

\maketitle

\begin{abstract}
We attempt to answer a question of D. C. Ravenel's:
``\textit{Infinite loopspace theorists, where are you?}'', or more prosaically:
 Can one start from the notion of Ring Spectrum Representing [Formal Group Laws] (the functor) and arrive at Complex Cobordism?
 Our answer is a tentative ``sort-of''.
\end{abstract}

\newtheorem{theorem}{Theorem}
\newtheorem{cor}[theorem]{Corrolary}
\newtheorem{lemma}[theorem]{Lemma}
\newtheorem{example}{Example}
\newcommand{\mg}{\mathcal{G}}
\newcommand{\Th}{\rlap{\rule{1pt}{0pt}\rule[.6ex]{.6em}{.4pt}}\mathrm{Th}}
\newcommand{\deloop}{\mathbf{B}}

\section{Introduction}
Study of complex cobordism goes back a long way, with the spectrum representing it being constructed by Pontrjagin and Thom.  Milnor showed that the cohomology ring $MU$ was {\em naturally} isomorphic to Lazard's ring $\Lambda$, while Landweber \cite{landweber} and Novikov \cite{novikov} show the stronger result that the Hopf algebra $(MU^*MU,MU^*)$ {\em naturally} represents Lazard's Hopf algebra.  Novikov further shows how to adapt a spectral sequence due to Adams, so that  complex cobordism gives a convenient tool for calculating many features of the Stable Homotopy Ring $\pi^S_*$; from another angle, there is an Atiyah-Hirzebruch spectral sequence $E^2_{*,*}\simeq H_p(\deloop U,\pi^S_q)\Rightarrow MU_{p+q}$; Brown+Peterson had already begun to show how to localize these techniques, and Chromatic Homotopy has more-or-less blossomed thence in many ways.  Anyone who wants to know the {\em real} story of these things should try Ravenel's excellent book \cite{ravenel}.

In the final section we make use of Betley's calculations \cite{betley}; his paper includes a passage I feel some need to echo:
\begin{quotation}
\noindent I believe that such a [.] result has been well known for many years but I could not find it in the literature.
\end{quotation}
Of course, this inability may be my own obtuseness, or Ravenel and I may be victims here of a mathematician's obscurantism the like of which I am myself too often guilty.  As a curious historical parallel, Steiner seems to have been the first one to even {\em pose} the isoperimetric problem, though it would certainly have been accessible to Euler, Lagrange, Leibniz, Newton\dots and Archimedes.

\subsection*{Notation} We write $\CP$ for infinite complex projective space, $\CP_k$ for its usual $2k$-skeleta, and will otherwise (try to) emphasize monoidal structures in superscripts: $\CP^{\times k}$ is a $k$-fold direct product, $\CP^{\smsh k}$ a $k$-fold smash, et.c.

We use $k!$ as the name of the symmetric group on $[k]$ (the von Neumann ordinal, for definiteness), whose cardinality is the usual significance of ``$k!$'';
Borel quotient is denoted $\bq $ and cofiber is $/$; in particular $*\bq G \simeq \deloop G$.
We will also use a generalized notion of Thom space, which {\em should} be familiar under the usual name of ``twisted half-smash'',
 the spaces approximating spectra representing twisted {\em homology}.
Curiously the only actual use of the phrase ``twisted half-smash'' I have yet seen in the liturature is singular also in that the twisting there ends up being trivial.
In any case, our twisted half smashes, or generalized Thom spaces, are denoted and defined as
$$ \Th(X,F) \simeq E/X $$ 
where $X \to E \to X$ is a section/retract with fibers $F\to E \to X$.
(The other obvious fibers are $\Omega F \to X \to E$).  
In particular, the Thom space of a group action in the pointed category is a cofiber of Borel quotients:
$$ \Th(\deloop G, X) \simeq (X\bq G) / (*\bq G)  .$$

\subsection*{Terminology and related}

We make frequent use of geometric arguments to study classifying spaces of {\em geometric things} and sometimes use well-known skeleta; nonetheless our perspective is {\em radically} homotopical:
 we do not work in any particular model categories
 (in particular, we never use notions of fibration or cofibration, nor mention weakness of equivalences),
 but in a setting where ``commuting diagrams'' come with (usually unnamed, but as-specific-as-sensible) homotopies, of whatever shape is required.  
In the metalanguage around Homtopy Type Theory, ``commutativity is data''.  
In fact what diagrams we use are, almost exclusively, coherent cubes, which are diagrams all of whose cubical facets are ``contractible'', to the extent that makes sense;
 and whenever we say ``colimit'' or ``cofiber'' or ``fiber'' we mean, exclusively, {\em homotopy}-colimit, {\em homotopy}-(co)fiber, et.c.

This does mean that such lovely things as differential cohomology and intersection cohomology just don't make sense
 (although cf. {\em cohesive} Homotopy Type Theory, as in \cite{schreiber} and used, e.g. in \cite{shulman});
 but we find this does not hamper the present investigation.

Lastly, within the usual Kelley category of CG[w]Hf spaces, one has Mather's Cube theorems; we call the {\em result} of Mather's theorems
 (and its finitary corollaries over any colimit) {\em distributivity}, after the obvious arithmetic analogue.  
Note that its use in homotopy is decades older than Mather's published proofs, underlying at least one classic construction of the Leray-Serre spectral sequence,
 (see also in particular the works of Tudor Ganea),
 and we use finitary distributivity as a general {\em axiom} of [unstable] homotopy.  
\footnote{Another echo rings in an aphorism of Rezk that distributivity is what promotes an $(\infty,1)$-category with distinct initial and terminal objects to {\em a Topos}, whatever controversy might otherwise surround that terminology}

\section{Algebraic motivation}
The place we start is: a Complex Oriented Commutative Ring Spectrum probably supports a Formal Group Law.  
That is, given a ring spectrum $\Lambda$ (suitable for Atiyah-Hirzebruch purposes) and a class 
 $g : \Lambda^2(\CP)$ extending the fundamental class $\eta: \Lambda^2(\sph^2)\sim \Lambda^0(*)\ni 1$,
 pull-back along tensor product of lines $\otimes: \CP\times\CP \to \CP$ induces a ring homomorphism
 $\Lambda^*(\CP) \to \Lambda^*(\CP\times\CP)$ which (due to the goodness of $g$ and collapsing spectral sequences) gives a formal group law
$$ \Lambda^*(\sph)[x] \to \Lambda^*(\sph)[x\otimes 1,1\otimes x] .$$
It may be worth remarking that the principal tools in the preceding are:
 the Atiyah-Hirzebruch spectral sequence collapses for an even-concentrated cohomology theory over an even-concentrated space;
 and the ordinary cohomology of $\CP$ is the ordinary polynomial ring $\ZZ[x]$;
 this feature it shares with the \textit{homology} of the odd-sphere loopspaces $\Omega \sph^{2n+1}$;
 indeed the Hopf Algebras $H^*(\CP)$ and $H_*(\Omega\sph^{2n+1})$ are isomorphic, which means it's hard to see
 (except that $\otimes$ is plainly special) why one should prefer {\em co}homology for FGL purposes.  
In other words, if this is already too much hypothesis in the direction of the conclusion for Ravenel's intent\dots \textit{H\'elas!}

{\em Having} chosen to reflect FGLs by complex-oriented ring spectra, one can clearly ask:
 what is the minimal complex-oriented ring spectrum?  
And this in turn means: what is the ring spectrum $\Lambda$ {\em generated by} a map $\CP \to \Lambda_2$?  
Answering such a question will proceed most simply if we choose a particular model (not to say ``model category''!) of ring spectra;
 and for this we choose the simplest available option, algebras of the $E_\infty$ operad in the homotopical monoidal category:
 ($E_\infty$ spaces, spectral smash product).

Using this model (cf. \cite{mayRings}) it is straightforward that our ring spectrum $\Lambda$ should be approximated by the spaces
$$ L_{2k} = \Th(\deloop k!,\CP^{\smsh k}) $$
with prespectrum maps induced by
$$ \sph^2 \smsh L_{2k} \to \CP \smsh L_{2k} \to L_{2k+2} $$
The remainder of this note is directed towards studying these $L_{2k}$ to show that they do in fact approximate, in the rough, {\em something}, and in the fine, compare it to our beloved $MU$; but, note: it is already looking difficult, in that there seems little reason for $L_{2k}$ to be homologically concentrated in even degree!

\section{Justifying the Presentation}

Perhaps it is reasonable to ask why the approximations $L_{2k}$ employ only a single operation; why there are no approximants $\Th(\deloop k!,\bigwedge_{i < k} L_{m(i)})$ in the rest of this study.  But the reason is simple: the $L_{2k}$ are already universal.  More: it is required, for any commutative ring spectrum $E$ and orientation $\CP\to E_2$, that we have {\em ring maps} $\Lambda_k \to E_k$; and ``ring maps'' {\em means} commutativity
$$ \begin{xy}\xymatrix{
\Th(\deloop k!, \Lambda_l^{\smsh k}) \ar[r] \ar[d] & \Th(\deloop k!,E_l^{\smsh k}) \ar[d]\\
\Lambda_{kl} \ar[r] & E_{kl}
}\end{xy} $$
Thus on the one hand, it is {\em necessary} that
$$\begin{xy} \xymatrix{
\Th(\deloop k!,\CP^{\smsh k}) \ar[r] \ar[d] & \Th(\deloop k!, E_2^{\smsh k}) \ar[d] \\
L_{2k} \ar[r] & E_{2k} } \end{xy} $$ commutes; and on the other hand, there is {\em exactly} one way for this to happen when the left arrow is an identity.

\section{Pursuing the product}
Since direct sum is plainly a useful $E_\infty$ construction on complex vector spaces (this is, after Segal, what makes $\deloop U$ an infinite loopspace), it is the most natural thing in geometry to consider the maps
$$ \bigoplus : \CP^{\times k} \bq  k! \to \deloop U(k).$$
We connect the Thom space $\Th(\deloop k!,\CP^{\smsh k})$ via this natural construction to the classic cobordism spectrum $MU$ using a fact, which surely is well-known but not often mentioned: the ordinary Thom spaces $MU(k)$ are the cofibers of the usual maps $\deloop U(k-1)\to \deloop U(k)$.  That is: the Thom space is described as the quotient (cofiber) of a (vector space or disk) bundle by the inclusion of the bundle of spheres (at-infinity or unit spheres); but for homotopy purposes, the important feature is the bundle of spheres itself, which is to say the map
$$ E \to \deloop U(k) $$ whose fibers are the unit spheres of the tautological bundle.  
The usual discussion of disk bundle and sphere-at-infinity acheives the necessary trivial-fibration$\circ$cofibration replacement in any of the classic model-category structures; however, the result of the construction is the (we now don't need to {\em say} ``homotopy'') cofiber.  Therefore, for the purposes of homotopy, it is sufficient to identify the type of the map $E\to U(k)$; but this is easy: the unit sphere is the (Borel, or ordinary) quotient $U(k)\bq U(k-1)$; the canonical action of $U(k)$ on this has Borel quotient equivalent to $\deloop U(k-1)$ and by functoriality of $\deloop$ the quotient {\em map} is the fiber map
$$ \sph^{2k-1} \to \deloop U(k-1) \to \deloop U(k) $$
Hence
\begin{lemma} $$ MU(k) \simeq \deloop U(k) / \deloop U(k-1) $$ \end{lemma}
From the source of the direct-sum map, we have another natural map
$$ \CP^{\times k}\bq k! \to \Th(\deloop k!,\CP^{\smsh k}) $$
which is, like $\deloop U(k)\to MU(k)$, a cofiber map.
More explicitly:
\begin{lemma} There is a pasting diagram of pushouts
$$\begin{xy}\xymatrix{ W^k \CP \bq k! \ar[r] \ar[d] & \CP^{\times k}\bq k! \ar[d] \\
\deloop k! \ar[r] \ar[d] & \CP^{\smsh k}\bq k! \ar[d] \\
{*} \ar[r] & \Th(\deloop k!, \CP^{\smsh k}) }
\end{xy} $$
\end{lemma} Specifically, $W^k \CP$ is the so-called ``fat wedge'', which we shall study more closely in the next sections.  But to summarize the state of things now, we have the following diagram, of which the verticals are cofiber sequences
$$\begin{xy}\xymatrix{
W^k \CP\bq k! \ar[d] & \deloop U(k-1) \ar[d] \\
\CP^{\times k} \bq k! \ar[r] \ar[d] & \deloop U(k) \ar[d] \\
\Th(\deloop k!,\CP^{\smsh k}) & MU(k) 
}\end{xy}$$
The next immediate objective is to construct a horizonal map between the cofibers by constructing a {\em lifting} along the top.
\section{The Cube of Coordinate Planes}
In any category with products, for any pointed object $* \to X$ and any $k:\NN$ there is a commuting $k$-cube of coordinate planes: the cube of pullbacks of $k$ inclusions
$$ X^{k-1} \to X^k $$ 
which themselves are the fibers {\em at} $*\to X$ of the $k$ defining structure maps $X^k \to X$.
(All of this is a belaboured way of saying: products should be as functorial as possible.)

Small instances of the coordinate planes cube are
$$ * \to X \qquad (k=1) 
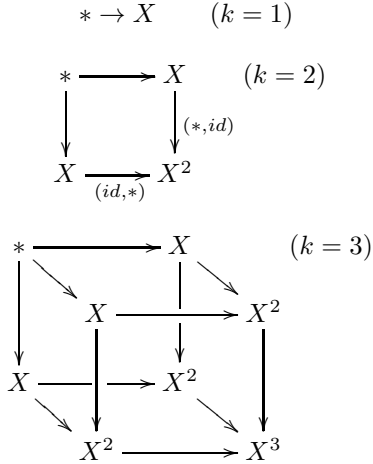
$$
$$ \begin{xy}\xymatrix{ {*} \ar[r] \ar[d] & X \ar[d]^{(*,id)} \\
X \ar[r]_{(id,*)} & X^2 } \end{xy} (k=2) $$
$$ \begin{xy}\xymatrix@=1.2em{
 {*} \ar[rr] \ar[dr] \ar[dd] &  & X \ar[dr] \ar[dd]|\hole\\
    & X \ar[dd] \ar[rr] &  & X^2 \ar[dd] \\
 X \ar[rr]|\hole \ar[dr]   &  & X^2 \ar[dr]  \\ 
  & X^2 \ar[rr] &  & X^3 }\end{xy} (k=3) $$
Though it is more generality than we need, different truncations of these diagrams (\dots the truncations should have a snappy name; I will call them ``blooms'') have colimits that generalize wedge sums, including, in particular, the ``fat wedge''; for the sake of convenience in notating $W^k X \bq k!$, we have chosen to index the fat wedges by the dimension of ambient cube defining them, rather than the degree of their petals.  This may be in conflict with other uses in the literature; please indulge me for the length of this note.
\begin{figure}[!h]
$$ \begin{xy}\xymatrix@=1.2em{
 {*} \ar[rr] \ar[dr] \ar[dd] &  & X \ar@{..>}[dr] \ar@{..>}[dd]|\hole\\
    & X \ar@{..>}[dd] \ar@{..>}[rr] &  & X\wsm X \ar@{..>}[dd] \\
 X \ar@{..>}[rr]|\hole \ar@{..>}[dr]   &  & X\wsm X \ar@{..>}[dr]  \\ 
  & X\wsm X \ar@{..>}[rr] &  & X\wsm X\wsm X }\end{xy} $$
\caption{\small Example: Wedge of three $X$s; the bloom looks like a sparse parsley inflorescence}
\end{figure}
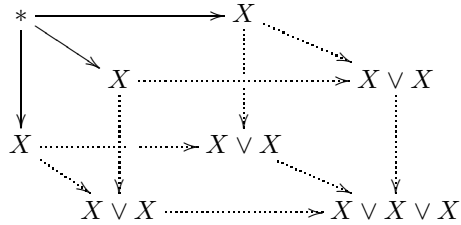
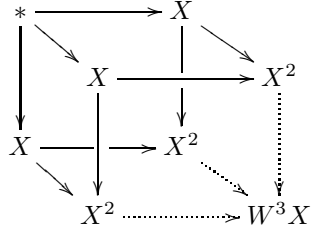
\begin{figure}[!h]
$$ \begin{xy}\xymatrix@=1.2em{
 {*} \ar[rr] \ar[dr] \ar[dd] &  & X \ar[dr] \ar[dd]|\hole\\
    & X \ar[dd] \ar[rr] &  & X^2 \ar@{..>}[dd] \\
 X \ar[rr]|\hole \ar[dr]   &  & X^2 \ar@{..>}[dr]  \\ 
  & X^2 \ar@{..>}[rr] &  & W^3 X }\end{xy} $$
\caption{\small Example: $W^3 X$; this bloom looks more like a tulip; in both examples, the dotted portions are generated as an initial Kan extension, to the cube, of the (solid-arrow) bloom}
\end{figure}

Our purpose in reminding the reader of this construction (who likely enough knows it already) is to prepare for the construction of the actual lifting, a map $W^k \CP\bq k! \to \deloop U(k-1)$; since $W^k\CP$ is a {\em colimit}, constructing a map out of it is accomplished by constructing a coherent cube, such as
$$ 
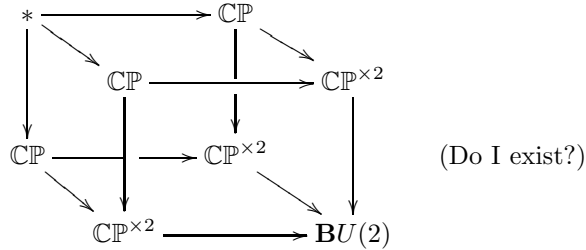
\begin{xy}\xymatrix@=1.2em{
 {*} \ar[rr] \ar[dr] \ar[dd] &  & \CP \ar[dr] \ar[dd]|\hole\\
    & \CP \ar[dd] \ar[rr] &  & \CP^{\times 2} \ar[dd] \\
 \CP \ar[rr]|\hole \ar[dr]   &  & \CP^{\times 2} \ar[dr]  & & \mbox{(Do I exist?)} \\ 
  & \CP^{\times 2} \ar[rr] &  & \deloop U(2) }\end{xy} $$
Now, this in itself is {\em not sufficient} for that we want a map on the Borel quotient $W^k\CP\bq k!$, of which $W^k\CP$ is just one fiber over $* \to \deloop k!$; and so we shall also have to take care that our constructed cube is suitably equivariant.

There isn't much choice about the maps towards the terminal vertex in the wanted diagram: we have direct sum, $\CP^{k-1}\to \deloop U(k-1)$, and anything hugely different is likely to be too forgetful.  We have already described the map $\deloop U(k-1)\to \deloop U(k)$ as the tautological sphere bundle over $\deloop U(k)$
$$ \sph^{2k-1} \to \deloop U(k-1) \to \deloop U(k) $$
which is to say that a point of $\deloop U(k-1)$ over $V:\deloop U(k)$ corresponds to {\em a unit vector in $V$}, which also articulates the sense in which the map $\deloop U(k-1)\to \deloop U(k)$ is direct sum {\em with $\CC$} and not just {\em any line}; in terms of the space $\deloop U(k-1)$ as classifier of $k-1$-dimensional complex vector spaces, a unit vector in $V$ corresponds to {\em the orthogonal complement} of that vector; and so a composite $\deloop U(k-2) \to \deloop U(k)$ further corresponds to {\em an ordered pair} of orthogonal unit vectors in $V$, and classifies the orthogonal complement of their span, et.c.  What matters for us is: the action of $k!$ on the relevant cubes {\em permutes these orders}; or, for now passing over the action of $k!$, the possible obstruction to homotopy commutativity in 
$$\begin{xy}\xymatrix{ \CP^{\times k-2} \ar[r] \ar[d] & \CP^{\times k-1} \ar[d] \\
\CP^{\times k-1} \ar[r] & \deloop U(k-1) 
}\end{xy}$$
is that the two routes have chosen orthogonal vectors in $V$, the unambiguous vector space specified by the unambiguous composite $\CP^{\times k-2} \to \CP^{\times k} \to \deloop U(k)$.  To specify the commutativity of the square, it is sufficient to specify a homotopy in unit vectors between these orthogonal two.  Happily, a basis for a complex space induces an inclusion of a real space with a basis, and in a real space with a basis, the positive orthant of the unit sphere does everything we need!  Being described in such dimension-agnostic terms, it is also clear how to generalize the argument to other layers of the cube, and commutativity with the already-known maps to $\deloop U(k)$.  And, hopefully, the equivariance of the construction is also sufficiently evident.

\begin{lemma}\label{liftingsquare} We therefore claim a commutative square
$$ \xymatrix{W^k\CP \bq k! \ar[r]\ar[d] & \deloop U(k-1)\ar[d] \\
\CP^{\times k} \bq k! \ar[r] & \deloop U(k)} $$ and hence a map 
$$ \Th(\deloop k!,\CP^{\smsh k}) \to MU(k) .$$
\end{lemma}

\section{Gluey Wedge, Take II}

The fiber of including a fat wedge in a product $W^k X \to X^k$, (and therefore the fiber of the map $W^k \CP \bq k! \to \CP^{\times k} \bq k!$ as well) is, by Distributivity, the $k$-fold self-join of the loopspace $(\Omega X)^{\star k}$; and so we have the small surprise: a fiber sequence
$$ \sph^{2k-1} \to W^k\CP \bq k! \to \CP^{\times k} \bq k! $$ 
In fact, we will argue that the square we already constructed by an interpolation argument (resulting in Lemma \ref{liftingsquare}), is a pullback.  We will work on the $k!$-space rather than the borel quotient, for now, but the equivariance of the argument should be clear.  To keep a handle on things, we'll start by naming the pullback itself
$$ \begin{xy}\xymatrix{
R_k \ar[r] \ar[d] & \deloop U(k-1) \ar[d] \\
\CP^{\times k} \ar[r]_\oplus & \deloop U(k)
}\end{xy} $$

To begin with, note once more that a point of $R_k$ over $x : \CP^{\times k} $ is specified by a unit vector {\em in} $\oplus x$.  Note furthermore that $x$ {\em is} (or might as well be) a {\em split} complex vector space: it comes with a prefered list of independent lines; this in turn gives the unit sphere in $\oplus x$, {\em canonically}, the structure of a join of (unpointed) circles:
$$ \begin{xy}\xymatrix{
{\underset{j}{\bigstar}\,\sph(x_j)} \ar[r] \ar[d] & R_k \ar[r] \ar[d] & \deloop U(k-1) \ar[d] \\
{*} \ar[r]_{(x_j|j=1\dots k)} & \CP^{\times k} \ar[r] & \deloop U(k) 
}\end{xy}$$
Moreover, this join itself is the colimit of the penultimate bloom within the product cube of the $k$ trivial maps $\sph(x_j) \to *$.  
Integrating this fiberwise decomposition over the base space $\CP^{\times k}$,
 the pullback $R_k$ is similarly the colimit of the bloom of various maps of total spaces;
 but this bloom lives in a cube again seen (with a little care) to be equivalently the cube of coordinate planes in $\CP^{\times k}$.  Consequently,
$$ R_k \simeq W^k \CP .$$

\section{Comparing with Cobordism}
\label{compar}
We should like to control the fibers or cofibers of $\Th(\deloop k!,\CP^{\smsh k}) \to MU(k)$.  The cofiber --- let us call it $Y_k$ --- is, of course, the Total Cofiber of a pullback square, which means it is also the cofiber of the comparison map between the base $\deloop U(k)$ and the pushout of the pullback structure maps:
$$\begin{xy}\xymatrix{
W^k\CP \bq k! \ar[r] \ar[d] & \deloop U(k-1) \ar[d] \ar@/^/[ddr] \\
\CP^{\times k} \bq k! \ar[r] \ar@/_/[drr] & P_k \ar[dr] \\
& & \deloop U(k) }\end{xy} $$
$$ Y_k \simeq \deloop U(k) / P_k \simeq MU(k) / \Th(\deloop k!,\CP^{\smsh k} ) $$
Distributivity now tells us that the {\em fibers} of the new map presenting $Y_k$ are the join of the vertical and horizontal fibers:
$$ \begin{xy}\xymatrix{
\sph^{2k-1} \times F_k \ar[r] \ar[d] & \sph^{2k-1} \ar[d] \ar@/^/[ddr] \\
F_k \ar[r] \ar@/_/[drr] & \sph^{2k-1} \star F_k \ar[dr] \\
& & {*} & \rlap{$\mbox{fibers over }\deloop U(k)$}
} \end{xy} $$
And so the main thing now is to study $F_k$.

This is perhaps a good moment to mention what $\CP^{\times k}\bq k!$ is the classifying space {\em of}.  As plainly as possible, there are a split extension of groups
$$ 1\to \mathbb{T}^k \to G \to k!\to 1 $$ 
and various (conjugate) inclusions
$$ \mathbb{T}^k \to G \to U(k) $$ 
realizing a maximal torus in $U(k)$.  Concretely, one may take $G$ generated by the diagonal matrices and the permutation matrices.  Continuing,
$$ \CP^{\times k}\bq k! \simeq \deloop G $$
$$ F_k \simeq U(k) \bq G \simeq (U(k) \bq \mathbb{T}^k)\bq k! $$
When I began working towards this note, I wanted to tell you how very-connected $\sph^{2k-1}\star F_k$ must be, but, alas, it is not!  There is a simple induction argument that each (complex) flag manifold $V_k = U(k)\bq \mathbb{T}^k$ is simply connected: there are fiber sequences
$$ V_k \to V_{k+1} \to \CP_k $$
and the base case $V_2 \simeq \CP_1$; this means that $\pi_1 F_k \simeq k!$, so that $\pi_2 \Sigma F_k \simeq C_2 = \mathrm{ab}(k!) $.

\section{Rescuing the Dream}\label{rescue}

Up to the end of Section \ref{compar}, we mostly avoided actually calculating the (co)homology of anything in any theory, because the types or spaces themselves have had plenty to say.  Many readers, however, will already have concluded many things about the homology of the spaces $L_{2k}$.  In particular, 
$$ H^*(L_{2k}; \QQ) \simeq (H^*(\CP^{\smsh k};\QQ))^{k!} $$
That is, the rational cohomology of $L_{2k}$ is the $k!$-invariants of the [already-torsion-free] homology of $\CP^{\smsh k}$ (because for finite $G$, all extensions of finitely-generated $\QQ[G]$-modules are split).

A result that should be well-known due to Borel+Serre \cite{borSer} that in ordinary cohomology, $H^*(\deloop U(k),\ZZ) \to H^*(\CP^{\times k};\ZZ)$ is {\em exactly} the inclusion of the $k!$-invariants; it will follow from this, since we have an explicit-enough factorization of the maps
$$ \CP^{\smsh k} \to \Th(\deloop k!,\CP^{\smsh k}) \to MU(k), $$
and from the multiplicative underpinning of the Thom isomorphism, that our map $L_{2k} \to MU(k)$ is an equivalence {\em modulo torsion}.

We can therefore, perhaps, give Ravenel the more concrete answer: Complex Cobordism is the initial {\em torsion-free} complex-orientable ring spectrum.

\section{Stabilization; or, $L_*$ presents a finite-type spectrum}

Having dismissed our homological hangups, before discussing what the {\em maps} do, let's consider the spaces themselves,
$$ \deloop k! \to \CP^{\smsh k}\bq k! \to \Th(Bk!,\CP^{\smsh k}).$$ 
The Cartan-Leray spectral sequence tells us
$$ H_p( k! ; H_q(\CP^{\smsh k},\mathbb{F})) \Rightarrow H_{p+q} (\CP^{\smsh k}\bq  k!) $$ 
while the cofiber map trivializes a direct summand $H_*(k!)\hookrightarrow H_*(\CP^{\smsh k}\bq k!)$, the image of the retract $\deloop k!$ of $\CP^{\smsh k}\bq k!$.
The $H_*(\CP^{\smsh k})$ happen to be sums of (possibly trivial) tensor products of the tautological representation;
 the rational case of this sequence (and sufficient evidence of its collapse) has already been used in \mbox{section \ref{rescue}}.

Betley \cite{betley} proves \begin{lemma}[Betley]\label{betstable} For any functor $T \in Ab^\Gamma$ of finite degree $d$, the natural map $$H_i( (n-1)! ; T ([n-1])) \to H_i( n! ; T([n])) $$ is an isomorphism provided $2i + d < n$.
\end{lemma} Here, $\Gamma$ is the category of {\em partial} maps of finite sets; the {\em degree} of such an Abelian functor is defined there in terms of the joint kernels of forgetting: a functor $T$ is of degree at-most $d$ whenever, for all $d$-part partitions $X = \coprod_j X_j$, the natural maps
$$ T(X) \to \bigoplus_i T(\coprod_{j\neq i} X_j) $$ 
are all injective.

For stable-homotopy purposes, we are interested in the homology groups $H_{2k+j}(\Th(k!,\CP^{\smsh k}))$ for $j\lesssim J$; that is, we are concerned only with the $2k+2J$-skeleton of the smash.  Particularly, the $2k+2J$-skeleton of $\CP^{\smsh k}$ consists of $\binom{J+k}{J}$ cells (apart from the basepoint).  More: in $\CP^{\smsh k}$, nontrivial cells are products of $k$ nontrivial cells in $\CP$, of even dimensions, at least $2$; if their total dimension is not more than $2k+2J$, then at most $J$ of the factors are of dimension {\em greater} than $2$; that is, if $k > J$, all the cells in the $2k+2J$-skeleton of $\CP^{\smsh k}$ are contained in an image of $$ \CP^{\smsh l} \smsh \sph^2 \smsh \CP^{\smsh k-l-1} \hookrightarrow \CP^{\smsh k} $$
and an easy amplification of this argument to $J+1$ smash factors leads us to the conclusion:
\begin{lemma} $H_{2k+j}(\CP^{\smsh k})$, for $j < 2J$, extends along $[k] : \Gamma$ as a functor of degree at most $J$.\end{lemma}
Consequently $H_p(k!;H_{2k+j}(\CP^{\smsh k}))$ stabilizes for $2 p + J \leq k$.

This establishes that the pages $E^2_{*,*}$ of the Cartan-Leray spectral sequences are stable in the sense necessary for $L_{2k}$ to approximate a finite-type spectrum; concluding that $L_{2k}$ themselves stabilize reduces to arguing that the prespectrum maps $\Sigma^2 L_{2k} \to \CP\smsh L_{2k} \to L_{2(k+1)}$ {\em realize} the morphisms shown in \mbox{Lemma \ref{betstable}} to become isomorphisms; but this already follows from the usual calculation of $H_i(\CP^{\smsh k})$. 

We therefore claim:
\begin{theorem} The spaces $L_{2k}$ with the maps $\Sigma^2 L_{2k} \to L_{2(k+1)}$ approximate a finite-type commutative ring spectrum $\Lambda$; this is the initial complex-oriented commutative ring spectrum.\end{theorem}

\end{document}